\documentclass[12pt,oneside]{amsart}
\usepackage{latexsym}
\usepackage{amsmath,amssymb,amsthm}
\usepackage{amscd}
\usepackage[dvips]{graphics}
\newcommand{\rp}{\mathbb{RP}}
\newcommand{\re}{\mathbb{R}}
\newcommand{\sph}{\mathbb{S}}
\newcommand{\pgl}[1]{\mathbf{PGL}(#1)}

\newcommand{\gl}[1]{\mathbf{GL}(#1)}
\newcommand{\na}{\nabla}
\newcommand{\sfrac}[2]{{\textstyle \frac{#1}{#2}}}
\newcommand{\D}{\displaystyle}

\newtheorem{prop}{Proposition} 
\newtheorem{cor}[prop]{Corollary}
\newtheorem{thm}{Theorem}

\newtheorem{lem}[prop]{Lemma}

\theoremstyle{remark}
\newtheorem*{rem}{Remark}
\newtheorem*{ack}{Acknowledgements}
\newtheorem*{notea}{Note added in proof}

\begin{document}

\author{John C. Loftin}
\title{Riemannian Metrics on Locally Projectively Flat Manifolds}
\maketitle

\begin{abstract}
The expression $-\frac1u u_{ij}$ transforms as a symmetric $(0,2)$ tensor under projective coordinate changes of a domain in $R^n$ so long as $u$ transforms as a section of a certain line bundle.  
On a locally projectively flat manifold, the section $u$ can be regarded as a metric potential analogous to the local potential in K\"ahler geometry.  
Let $M$ be a compact locally projectively flat manifold. We prove that if $u$ is a negative section of the dual of the tautological bundle such that $-\frac1u u_{ij}$ is a Riemannian metric, then $M$ is projectively equivalent to a quotient of a bounded convex domain in $R^n$. The same is true for such manifolds $M$ with boundary if $u=0$ on the boundary. This theorem is an analog of a result of Schoen and Yau in locally conformally flat geometry.
The proof involves affine differential geometry techniques developed by Cheng and Yau. 
\end{abstract}

\section{Introduction}
A locally projectively flat manifold is a manifold built of coordinate charts in $\rp^n$ with transition maps in $\mathbf{PGL}(n+1,\mathbb{R})$.  In Thurston's language, this defines an $(X,G)$ manifold, with homogeneous space $X=\rp^n$ and Lie group $G=\mathbf{PGL}(n+1,\mathbb{R})$.  (From now on we drop the dependence on $\re$ from the Lie group notation.)  We will actually be concerned with a slightly smaller class $\mathcal{S}$ of manifolds $(X,G)$ manifolds for $X=\sph^n=(\re^{n+1}\setminus 0)/\re^+$ and $G=(\gl{n+1})/\re^+$. For a locally projectively flat manifold $M$, being able to lift the projective structure to a structure in $\mathcal{S}$ is equivalent to the existence of a tautological bundle $\tau$ over $M$ analogous to the tautological bundle over $\rp^n$ such that $\tau \to M$ is oriented (see Section \ref{lo-ni} below).

Now consider a negative section $u$ of the dual bundle $\tau^*$ over such a manifold $M$. If $u$ is strictly convex and satisfies the differential equation 
\begin{equation}
 \det (u_{ij})=\left(-\frac1u\right)^{n+2},\label{affine-sph-eq}
\end{equation}
 then the section $-\frac1u$ of the tautological bundle $\tau$ traces out a hyperbolic affine sphere in the total space of $\tau$. (Note: for a general hyperbolic affine sphere, we should replace this equation by $\det (u_{ij})=\left(\frac1{Lu}\right)^{n+2}$ for $L$ a negative constant.  We normalize by scaling along the fibers of $\tau$ to make $L=-1$.) For $M$ compact (without boundary), this implies that the universal cover of $M$ is projectively equivalent to a bounded convex domain \cite{loftin01}.  The reason is that if we look at the image of the developing map of $M$, $-\frac1u$ induces a complete immersed hyperbolic affine sphere in $\re^{n+1}$ with complete affine metric.  Then the fundamental work of Cheng and Yau \cite{cheng-yau86} shows that this affine sphere must be asymptotic to a cone over a bounded domain $\Omega$, and thus the universal cover of $M$ is projectively equivalent to $\Omega$. (The basic results on hyperbolic affine spheres are contained in Cheng-Yau \cite{cheng-yau77,cheng-yau86} and also Calabi-Nirenberg \cite{calabi-nirenberg74}. See also \cite{gigena81,li90,li92,li-simon-zhao,sasaki80}. Consequences of Cheng-Yau's work for locally projectively flat manifolds are found in \cite{loftin01} and also Labourie \cite{labourie97}.)

\begin{thm} \label{main-thm}
Let $M^n$ be a compact locally projectively flat manifold without boundary. The following statements are equivalent:
\begin{enumerate}
\item $M$ admits an oriented tautological bundle $\tau$ and there is a negative section $u$ of the the dual bundle $\tau^*$ whose Hessian matrix $u_{ij}>0$ (derivatives taken with respect to some $\rp^n$ coordinates).
\item $M^n$ is projectively equivalent to $\Omega/\Gamma$, where $\Omega \subset \subset \re^n \subset \rp^n$, and $\Gamma \subset \pgl{n+1}$ acts discretely and properly discontinously on $\Omega$.
\end{enumerate}
\end{thm}

If $M$ satisfies condition 2\ above, we say $M$ is {\it properly convex}.  In particular the developing map from the universal cover $\tilde{M}$ is injective and its image is the domain $\Omega$ which is properly contained in an affine $\re^n \subset \rp^n.$ We note that 2\ implies 1\ is a consequence of Cheng and Yau's work on Monge-Amp\`ere equations \cite{cheng-yau77}.  Our contribution is then 1\ implies 2.

If we have $u$ a negative strictly convex section of $\tau^*$, one approach to showing that $M$ is properly convex is to try to perturb $u$ to a section satisfying the affine sphere equation \ref{affine-sph-eq}. In other words, find a positive function $\phi$ such that $\phi u$ satisfies
\begin{equation}
 \det ((\phi u)_{ij})=\left(-\frac1{\phi u}\right)^{n+2}, \qquad (\phi u)_{ij}>0. \label{perturb-eq}
\end{equation}
Although we do not use this approach (in fact we can solve this equation as a corollary of the main theorem and Cheng-Yau \cite{cheng-yau77}), we remark that the idea of deforming an object to a canonical object by PDE techniques and then using the canonical object to obtain geometric or topological information has several parallels. For a complex manifold, there is the famous K\"ahler-Einstein equation, in which one perturbs a K\"ahler metric $g_{i\bar{\jmath}}$ by the complex Hessian $\phi_{i\bar{\jmath}}$ of a function $\phi$. 

Analogous statements for other $(X,G)$ manifolds are perhaps less well known.  $(X,G)$ manifolds are more rigidly constructed than complex manifolds; so it is therefore not surprising that quite precise statements on the nature of the universal cover of the given manifold can be obtained by perturbing a given metric to a more ``canonical'' metric. For locally conformally flat manifolds, Schoen and Yau \cite{schoen-yau88,schoen-yau94} prove that for such a compact manifold $M$, the existence of a locally conformally flat metric $g$ with scalar curvature $R=1$ implies that $M$ is a quotient of $\Omega \subset \sph^n$, where $S^n \setminus \Omega$ has small Hausdorff dimension.  Then using Schoen's solution to the Yamabe problem \cite{schoen84,schoen-yau94}, any conformally flat metric with $R>0$ has in its conformal class a metric of constant positive scalar curvature.  Thus the closed condition $R=1$ can be replaced by the open condition $R>0$ to conclude that $M$ is a quotient of some $\Omega$ as above.

Also, Cheng and Yau in \cite{cheng-yau82} use the analog of the K\"ahler-Einstein equation on affine flat manifolds to find various conditions in terms of a given Hessian metric to perturb it to the analog of a K\"ahler-Einstein metric.
Then they use properties of this metric to show that the universal cover is affinely equivalent, in various cases, either to all of $\re^n$ or to a pointed convex cone in $\re^n$.  We should also point out the result of Shima \cite{shima81} that any compact affine flat manifold which admits a Hessian metric is a quotient of a convex domain, and that Shima and Yagi \cite{shima-yagi97} have extended this result to affine flat manifolds with complete Hessian metrics. 

Instead of solving Equation \ref{perturb-eq} above, we prove our main result directly by using the affine differential geometry techniques developed by Cheng and Yau in \cite{cheng-yau86}.
This approach has the advantage of solving the related problem when $M$ is allowed to have boundary and Dirichlet condition $u=0$ on the boundary (Theorem \ref{thm-boundary} below).

\begin{rem} Another approach to this problem is to use Andrews's affine normal flow \cite{andrews96} to flow our initial hypersurface as above to an affine sphere.  In other words, the initial hypersurface as the radial graph of $-\frac1u$ should converge under the affine normal flow  to a homothetically expanding affine sphere.  Then Cheng-Yau's results as above will apply to show our manifold must be properly convex. Similar estimates to those in \cite{andrews96} should provide this result, but the case with Dirichlet boundary $u=0$ may cause difficulty.
\end{rem}

\begin{notea}
There is a simple relationship between the metrics we consider here and Hessian metrics on the total space of the tautological bundle.  This relationship allows another approach to the results here, and relates the affine normal flow of a hypersurface to a natural flow of Hessian metrics on the total space of $\tau$.  These results will appear in a future paper.
\end{notea}

\begin{rem} Recently, Choi \cite{choi99} has decomposed any locally projectively flat manifold $M$ into a number of geometric pieces with varying convexity properties.  In particular any such manifold can be split into a number of pieces with totally geodesic boundary, each piece either affine flat or at worst $(n-1)$-convex.  (Here 1-convex is the usual notion on convexity and, roughly speaking, $m$-convex means replacing the line segments in the usual definition of convexity with $m$-simplices. See \cite{choi99} for details.)
What we have proved here is that if there exists an appropriate section $u$ on $M$ which is strictly convex, then $M$ is a quotient of a bounded convex domain.  It is tempting to think that relaxing the convexity properties of $u$ will allow us to use more types of pieces in Choi's decomposition. For example, the existence of an affine flat piece may correspond to the Hessian of $u$ becoming degenerate and the presence of different $m$-convex pieces may signify that we should allow up to $m-1$ eigenvalues of $u_{ij}$ to become negative. 
\end{rem}

\begin{ack}
The author would like to thank Profs.\ S.-T.\ Yau, R.\ Schoen and W.\ Goldman for their interest in this work, and the referee for useful suggestions.
\end{ack}

\section{Loewner-Nirenberg's metrics} \label{lo-ni}
Loewner and Nirenberg in \cite{loewner-nirenberg} studied metrics of the form $-\frac1u u_{ij}$  on bounded convex domains, where $u$ is a convex negative function and $u_{ij}$ denotes the Hessian matrix. If we treat $u$ not as a function but as a section of a certain line bundle, then $-\frac1u u_{ij}$ transforms as a Riemannian metric under projective coordinate transformations. 

We introduce the tautological bundle to explain how $u$ transforms.
Consider the coordinates of our domain $\Omega$ to be inhomogenous projective coordinates $(t^1,\dots,t^n,1)$ in $\re^{n+1}$. Then a section $s$ of the tautological bundle $\tau$ is naturally represented by the radial graph $s(t^1,\dots,t^n,1) \subset \re^{n+1}$.  Let a matrix $A=(a_i^j)$ (normalized so that det$A=\pm1$) act on $\re^{n+1}$:
$$ \tilde{x}^j= a_i^jx^i, \quad 1\le i,j\le n+1$$
Then if $x= s(t^1,\dots,t^n,1)$, the formula 
$$\tilde{x}(x)=\tilde{s}(\tilde{t}^1,\dots,\tilde{t}^n,1)$$
defines both the coordinate transform to new homogeneous coordinates $\tilde{t}^j$ and the transition function for the tautological bundle $\tau$. 
In order to make the Hessian $u_{ij}$ transform correctly, we require $u$ to be a section of $\tau^*$ the dual of the tautological bundle.  Also Loewner and Nirenberg find that $(\det u_{ij})^\frac1{n+2}$ transforms as a section of the same bundle, and thus are led to study solutions to Equation \ref{affine-sph-eq}
$$\det (u_{ij})=\left(-\frac1u\right)^{n+2}$$
and the projectively invariant Riemannian metric $-\frac1uu_{ij}\,dt^idt^j$ its solution provides.
A deeper understanding of these objects is afforded by the affine differential geometry, which we outline below in Section \ref{aff-diff}.

But first we use Loewner-Nirenberg's analytic formulation to put a few necessary conditions on a locally projectively flat manifold $M$. Notice that in order for a metric of the form $-\frac1u u_{ij}$ to be defined on $M$, we must ensure that an appropriate analog of the tautological bundle exists on $M$.  Also, the requirement that $u$ be negative only makes sense if $\tau$ is oriented. 
Thus we must restrict our discussion to a slightly smaller class of manifolds, those which admit a tautological line bundle which is oriented. A more geometric reason for this requirement is that below we want to consider the geometry induced by the positive section $-\frac1u$ of $\tau$. As we see in Proposition \ref{convexity} below, the positivity of the Hessian $u_{ij}$ is equivalent to the graph of this section's being strictly convex as a hypersurface in the total space of $\tau$.

We briefly outline how to construct such a $\tau$.  Let $\mathbf{SL}^\pm(n+1)$ denote the group of all real $(n+1)\times(n+1)$ matrices with determinant $\pm1$. Under the natural projection, $\mathbf{SL}^\pm(n+1)$ is always a double cover of $\pgl{n+1}$. Blow up the origin in $\re^{n+1}$ to form the total space $\mathcal{T}$ of the tautological bundle over $\rp^n$. Lift the transition maps on $M$ from $\pgl{n+1}$ to $\mathbf{SL}^\pm(n+1)$ (if possible) to patch together open sets composed of lines in $\mathcal{T}$.  This forms the total space of a tautological bundle $\tau$ over $M$.  Finally we require that $\tau$ be oriented.  

(These two requirements together, the possibility of making a consistent choice of lifts and the orientability of $\tau$, are equivalent to being able to lift our $(\rp^n,\pgl{n+1})$ structure to an $(\sph^n,(\gl{n+1})/\re^+)$ structure, where $\sph^n$ is defined as the quotient space $(\re^{n+1}\setminus 0)/\re^+$.)

Also note that requiring $M$ to admit such a $\tau$ is no restriction in our case, as it is easy to show that any convex manifold admits an oriented tautological bundle. See e.g.\ \cite{loftin01} for details. 

\section{Affine differential geometry} \label{aff-diff}
In this section, we recall some basic formulas from affine differential geometry which will be used in the final section to complete the proof.  A basic reference for this is Nomizu-Sasaki \cite{nomizu-sasaki}. To any strictly convex hypersurface $H$ in $\re^{n+1}$ we can associate a transverse vector field $\xi$ which is invariant under unimodular affine transformations. In other words, for $\phi(x)=Ax+b$, $\left|\det A\right|=1$, then $\xi_{\phi(y)}=\phi_*(\xi_y)$ for any $y \in H$.  If the affine normal $\xi$ equals the position vector $x$, then $H$ is a hyperbolic affine sphere.  Furthermore, if $H$ is the radial graph of $-\frac1u$, then the equation for $\xi=x$ in terms of $u$ is Equation \ref{affine-sph-eq} above.

The affine normal, via the affine structure equations, provides us with tensors on $H$ which are also invariant under unimodular affine transformations.  The affine metric is a Riemannian metric, the cubic form $A$ is a symmetric $(0,3)$ tensor, and the affine shape operator $B$ is a section of $T^*\otimes T$.  In our formulation above, the section $-\frac1u$ embeds $M$ locally as a hypersurface in $\re^{n+1}$ with transition maps in $\mathbf{SL}^\pm(n+1)$.  Therefore, by the invariance of the affine normal, the section $u$ endows $M$ with the affine metric and the tensors $A$ and $B$.

Below we will need a formula for the affine metric in local coordinates. 
If we locally represent $H$ by the Euclidean graph of a function $f$, i.e.\ $H=\{(x^1,\dots,x^n,f(x))\}$, then the metric is given by
\begin{equation} \label{defaffinemet}
(\det f_{k\ell})^{-\frac1{n+2}}f_{ij} \, dx^i dx^j, \quad f_{ij}=\frac{\partial^2f}{\partial x^i \partial x^j}.
\end{equation}

The {\it affine conormal vector} is a vector in the dual space $\re_{n+1}$ defined at a point $x \in H$ by 
\begin{eqnarray*}
\nu_x(\xi)&=&1, \\
\nu_x(X)&=&0, \quad X \in T_x(H).
\end{eqnarray*}
For a linear map $A\in \mathbf{SL}^\pm(n+1)$, consider $\nu_{Ax}$ the conormal map at the point $Ax$ in the hypersurface $A(H)$. Then $\nu_{Ax}=(A^\top)^{-1}\nu_x$.   Therefore, if we again view $M$ locally as a hypersurface via the section $-\frac1u$, $\nu$ is well-defined on $M$ as a map to a section of the tautological bundle of the projective dual manifold of $M$.

We will need below a formula for the second covariant derivative of $\nu$. Choose normal coordinates with respect to the affine metric at a point $x\in H$.  Then we have at that point (see e.g.\ \cite[p.\ 66]{nomizu-sasaki})
$$ \nu_{,ij}=-\sum_k A_{ijk} \nu_{,k} - B_{ij}\nu.$$
Here $A_{ijk}$ is the cubic form, $B_{ij}$ is the shape operator with an index lowered by the metric, and the comma denotes covariant differentiation.


There is another natural conormal map, that with respect to the position vector. In this case we assume the position vector $x$ is transverse to $H$ for $x\in H$. We define the {\it centroaffine conormal vector} $\mu$ by 
\begin{eqnarray*}
\mu_x(x)&=&1, \\
\mu_x(X)&=&0, \quad X \in T_x(H).
\end{eqnarray*}
As above for $\nu$, $\mu$ is well defined on $M$ and provides a map to a section of the tautological bundle of the projective dual manifold of $M$.

It is also useful to have representations of $\nu$ and $\mu$ in local coordinates. Again write $H=\{(x^1,\dots,x^n,f(x))\}$. Then in the dual coordinates in $\re_{n+1}$,
\begin{eqnarray*}
\nu&=&\det(f_{ij})^{-\frac1{n+2}}(-f_1,\dots,-f_n,1), \\
\mu&=&-\frac1{f_jx^j-f}(-f_1,\dots,-f_n,1).
\end{eqnarray*}
For $f$ strictly convex, $v=f_jx^j-f$ is called the {\it Legendre transform} of $f$ and is thought of as a function in the variables $y_j=f_j$.  Then we may consider $(-y_1,...,-y_n,1)$ to represent inhomogeneous projective coordinates on the projectivization ${\mathbb P}(\re_{n+1})$.  In this coordinate patch, $-\frac1v$ and $\det(f_{ij})^{-\frac1{n+2}}$ locally represent sections of the tautological bundle on the projective dual manifold, and in particular, are sections of the same line bundle over $M$.  Also note that $\det(f_{ij})^{-\frac1{n+2}}$ is the last component of the vector $\nu$.

Similar reasoning shows why $-\frac1u$ and $\det(u_{ij})^{\frac1{n+2}}$ are sections of the same bundle.  If we consider locally the hypersurface $H$ given by the radial graph of $-\frac1u$, then $H=\mu(\mu(H))$ is of course the radial graph of $-\frac1u$, while $\nu(\mu(H))$ is the radial graph of $\det(u_{ij})^{\frac1{n+2}}$.


Finally we consider centroaffine geometry. In other words, we take structure equations with respect to the position vector, not the affine normal.  As long as the position vector of $x$ is transverse to $H$ for each $x \in H$, we can write
\begin{equation}
D_X\,x_*(Y)=x_*(\nabla^c_X Y) + h(X,Y)x, \quad X,Y \in T_x(H). \label{centroaffine}
\end{equation}
Here $D$ is the usual flat connection on $\re^{n+1}$, $\nabla^c$ is a connection on the tangent bundle $T(H)$, and $h$ is a symmetric $(0,2)$ tensor.  Again, in our situation, the centroaffine connection $\nabla^c$ descends to $M$.  The geodesics with respect to $\nabla^c$ are precisely the geodesics with respect to the $\rp^n$ structure on $M$, i.e.\ straight lines in any projective coordinate chart.

\section{The Legendre transform} \label{legendre-sect}

The estimates below depend heavily on properties of the Legendre transform.  We explain the basic calculations we will need for the reader's convenience.  Recall that the Legendre transform of a strictly convex function $f(x)$ is defined by 
\begin{equation} 
v+f=\frac{\partial f}{\partial x^i}\,x^i. \label{def-legendre}
\end{equation}
Here $v$ may be considered a function both of the variables $y_i=\frac{\partial f}{\partial x^i}$ and also of the variables $x^i$. Note that the strict convexity of $f$ insures that the map $x\mapsto y$ defined by $y_i=\frac{\partial f}{\partial x^i}$ sends any convex domain in $x$ diffeomorphically to a convex domain in $y$.  $v$ is considered primarily to be a function of the $y_i$.

Apply $\frac{\partial}{\partial x^j}$ to (\ref{def-legendre}) to show 
\begin{equation} \label{gradv}
\frac{\partial v}{\partial x^j}=x^i\frac{\partial^2f}{\partial x^i\partial x^j}.
\end{equation}
Also, we will need this lemma
\begin{lem} \label{positive}
If $f$ is strictly convex and has a minimum at the origin, then $v+f \ge 0$ on any domain in the $x$ variables star-shaped around the origin.
\end{lem}

\begin{proof}
It is clear from the definition (\ref{def-legendre}) that $v+f=0$ at the origin.  By changing coordinates it is enough to prove that $v+f\ge0$ at $q=(1,0,\cdots,0)$, assuming that $q$ is in our star-shaped domain.  Simply compute 
$$\frac{\partial}{\partial x^1} (v+f)= x^k\,\frac{\partial^2f}{\partial x^k \partial x^1} + \frac{\partial f}{\partial x^1}$$
and note that this is positive for all points on the line segment between the origin and $q$ since $f$ is strictly convex.
\end{proof}

\section{Convexity of a radial graph}
In this section we prove a proposition showing that if $u$ is strictly convex, then the radial graph of $-\frac1u$ is strictly convex as a hypersurface in $\re^{n+1}$.  We will be applying affine differential geometry techniques to this hypersurface, and we must check it is strictly convex to make sure the affine metric is positive definite.

\begin{prop} \label{convexity}
Let $u$ be a negative function defined on a domain $\Omega \subset \re^n$.  The following are equivalent.
\begin{enumerate}
\item $u$ is a strictly convex function on $\Omega$.
\item The radial graph of $-\frac1u$ is strictly convex as a hypersurface in $\re^{n+1}$ and the hypersurface points away from the origin (i.e.\ any tangent hyperplane puts the origin in the opposite half-space from a neighborhood in the hypersurface near the point of tangency).
\end{enumerate}
\end{prop}

\begin{proof}
We consider the centroaffine second fundamental form given by $h$ in Equation \ref{centroaffine}.  It is easy to check that $h$ is a nonzero scalar multiple of the usual second fundamental form [as is the affine metric or any other tensor defined by a transverse vector field in as in (\ref{centroaffine})].
Therefore, we only check that $h$ is definite. It is easy to check that the hypersurface points away from the origin.

Let $x$ be the position vector
$$ x=-\frac1u(t^1,\dots,t^n,1).$$
Then the tangent space of the image of $x$ is spanned by 
$$ \frac{\partial x}{\partial t^i} = x_* \! \left(\frac\partial{\partial t^i}\right)= -\frac1u \frac{\partial u}{\partial t^i}\, x + \left(0,\dots,0,-\frac1u,0,\dots,0\right),$$
where the only nonzero slot in the last term is the $i^{\rm th}$ one. Then 
\begin{eqnarray*}
&&D_{\frac\partial{\partial t^j}}\, x_* \! \left(\frac\partial{\partial t^i}\right)
= \frac\partial{\partial t^j}\,x_*\! \left(\frac\partial{\partial t^i}\right) \\
&&= \frac1{u^2} \frac{\partial u}{\partial t^j} \frac{\partial u}{\partial t^i} x
-\frac1u\frac{\partial^2 u}{\partial t^i \partial t^j} x - 
\frac1u \frac{\partial u}{\partial t^i} \frac{\partial x}{\partial t^j}
+ \left(0,\dots,\frac1{u^2} \frac{\partial u}{\partial t^j},0,\dots,0\right).
\end{eqnarray*}
The part of this expression in the span of $x$ is just
\begin{equation}
h\left(\frac{\partial}{\partial t^j}, \frac\partial{\partial t^i}\right)=-\frac1u\frac{\partial^2 u}{\partial t^i \partial t^j}\,\, , \label{c-metric}
\end{equation}
which is clearly positive definite.
\end{proof}

\begin{rem}
A more involved computation shows that the affine metric of the radial graph of $-\frac1u$ is given in the $t$ coordinates by
\begin{equation}
 \left( \det \frac{\partial^2 u}{\partial t^\ell \partial t^k}\right)^{-\frac1{n+2}} \, \frac1{u^2}\, \frac{\partial^2 u}{\partial t^i \partial t^j}\,dt^i\,dt^j. \label{a-metric}
\end{equation}
Also, Calabi showed in \cite{calabi72} that the affine metric of the (Euclidean) graph of the Legendre transform of $u$ is given by
\begin{equation}
 \left( \det \frac{\partial^2 u}{\partial t^\ell \partial t^k}\right)^{\frac1{n+2}} \, \frac{\partial^2 u}{\partial t^i \partial t^j}\,dt^i\,dt^j. \label{l-metric}
\end{equation}
Therefore if any two of these three metrics (\ref{c-metric}), (\ref{a-metric}), (\ref{l-metric}) are naturally isometric, then all three are, and $u$ must satisfy  Equation \ref{affine-sph-eq}. In this case, both the graphs in question are affine spheres (by the formulation in \cite{calabi72} and a duality result in Gigena \cite{gigena81}).
\end{rem}

\section{The main theorem}
In Cheng and Yau's work on affine spheres, a key point is that many of the various affine invariants, when contracted with the metric, are bounded on a complete affine sphere.  For example, Calabi's estimates \cite{calabi72,cheng-yau86} show the cubic form must be bounded whenever the affine metric is complete.  In the present case, the fact we are on a compact manifold bounds all of the affine invariants automatically.
The particular gradient estimate below follows the treatment in Li \cite{li92}.

Consider the section $u$ lifted to the universal cover $\tilde{M}$ of $M$. The the radial graph of the section $-\frac1u$ of $\tau$ is, by use of the developing map, an immersed locally strictly convex hypersurface in $\re^{n+1}$.  
Choose a point $p$ in $\tilde{M}$ and an inhomogeneous coordinate patch 
$\mathcal{O}$ containing $p$ so that $p$ is the origin and $u$ has a local minimum at~$p$. For this last statement, note that the choice of $\mathcal{O}$ also provides us with a local frame for the line bundle $\tau^*$.

Define $\Omega_h$ to be the connected component of the set 
$$\{t\in \mathcal{O}\cap \tilde{M}:\ u(t)<-\frac1h\}$$ 
containing the origin $p$.
We will show that for each $h>0$, $\Omega_h$ is properly contained in $\tilde{M}$ and that they exhaust $\tilde M$ as $h\to \infty$.
Since $u$ is convex and negative in the coordinate patch $\mathcal{O}$, and has its minimum at $p=0$, each $\Omega_h$ must be bounded in $\mathcal{O}$ and strictly convex.  Therefore, $\tilde{M} \subset \mathcal{O}$ must be convex, and since it supports a convex negative function $u$ with minimum at the origin, this inclusion is proper.

\begin{lem} \label{maxh-lemma}
Assume $\Omega_h$ is not properly contained in $\tilde{M}$ for some $h<\infty$. Then there is a minimal value of $h$ such that $\Omega_h$ is not properly contained in $\tilde{M}$.
\end{lem}
\begin{proof}
It is obvious that $\Omega_h \subset \subset \tilde{M}$ for $h$ near $-\frac1{u(p)}$.  We need to show this condition is open. Consider $h$ such that $\Omega_h \subset \subset \tilde{M}$. For a point $p'$ in $\partial \Omega_h$, the ray in $\mathcal{O}$ from $p$ to $p'$ is a geodesic for the real projective structure, and of the centroaffine connection $\nabla^c$.  Since $\nabla^c$ descends to $M$ and $\partial \Omega_h$ is compact, we can extend each ray a small distance.  Then since $u$ is increasing along each ray, we can extend $\Omega_h$ to $\Omega_{h+\epsilon}$ for small $\epsilon$.
\end{proof}

We assume such a minimal $h$ exists.  It is clear that $\Omega_h \subset \subset \mathcal{O}$ and $\Omega_h \subset \tilde{M}$, but there is some point $y \in \partial \Omega_h \subset \mathcal{O}$ with $y$ not in the closure of $\Omega_h$ in $\tilde{M}$. 

Since $\Omega_h$ is convex, we can define $f$ to be the function whose Euclidean graph is the radial graph of $-\frac1u$. 
\begin{equation}
(x^1,\cdots,x^n,f(x))=(-\frac{t^1}{u(t)},\cdots,-\frac{t^n}{u(t)},-\frac1{u(t)}) \quad \mbox{for} \quad t \in \Omega,
\end{equation}
Notice that $f<h$ on $\Omega_h$ and $f\to h$ at $\partial \Omega_h$. The $x^i$ are the coordinates for the Euclidean graph $f$ and $\nabla$ is the Levi-Civita connection with respect to the affine metric. We use the following gradient estimate to prove a contradiction.

\begin{prop}
Let $v$ be the Legendre transform of $f$. Then
$$\frac{\|\na v\|}{-v} \le C(h-f)^{-\frac12} \quad \mbox{on}\quad \Omega_h,$$
where $C$ depends only on $h$ and the geometry of $M$ and the section $u$.
\end{prop}
\begin{proof}
It is convenient to introduce the function $w=-\frac1v$. Of course $\frac{\na w}{w}=\frac{\na v}{-v}$.

Let $r$ be the geodesic distance function from $p$ with respect to the affine metric on the graph of $f$. 
Consider the following $\phi$ defined on $\Omega_h \cap B(p,a)$:
$$ \phi = (a^2-r^2)(h-f)^{\frac12}\frac{\|\na w\|}{w}. $$
(Here $a$ is a positive constant, $B(p,a)$ is the geodesic ball of radius $a$ centered at $p$, and the norm and gradient are with respect to the affine metric.)

At the point $q$ where $\phi$ achieves its maximum (which must be nonzero), $\na \phi=0$. (We can take care of the case when $q$ is on the cut locus, and therefore $\phi$ is not smooth, by Calabi's technique \cite{calabi57}.) Choose at $q$ normal coordinates with respect to the affine metric so that at $q$
\begin{eqnarray*}
w_{,1} &=& \|\na w\|>0, \\
w_{,i} &=& 0 \quad \mbox{for} \quad i\ge 2,
\end{eqnarray*}
In these coordinates we compute $\phi_{,1}=0$ at $q$:
\begin{equation} \label{gradphi}
\begin{array}{l}
 \D \sfrac12(a^2-r^2)(h-f)^{-\frac12}\frac{f_{,1}w_{,1}}{w} 
+ (r^2)_{,1}(h-f)^{\frac12}\frac{w_{,1}}{w} \\
\D {}- (a^2-r^2)(h-f)^{\frac12}\frac{w_{,11}}{w} 
- (a^2-r^2)(h-f)^{\frac12}\left( \frac{w_{,1}}{w} \right)^2 =0
\end{array}
\end{equation}
We can relate the fourth term to $\phi^2$: this term will be used to do the estimate. The first term, though potentially bad since it blows up, will fortunately be nonnegative by using properties of the Legendre transform (Section \ref{legendre-sect} above).  For the second term, note that 
\begin{equation} \label{second-term}
(r^2)_{,1}\le 2a
\end{equation}
 since $r \le a$. Analyzing the third term is the key estimate.

To take care of the first term in (\ref{gradphi}), let $(f^{ij})$ denote the inverse matrix of $(f_{ij})=(\frac{\partial^2f}{\partial x^i \partial x^j})$ and compute
\begin{equation} \label{first-term}
\begin{array}{l}
\D \frac{f_{,1}w_{,1}}w = \frac{f_{,1}v_{,1}}{-v}=\frac1{-v}\sum_i f_{,i}v_{,i}
\\ \D = \frac1{-v}\sum_{i,j} (\det f_{kl})^\frac1{n+2}f^{ij} \, \frac{\partial f}{\partial x^i}\,\frac{\partial v}{\partial x^j} \\
\D =\frac1{-v} (\det f_{kl})^\frac1{n+2} \sum_i x^i \frac{\partial f}{\partial x^i}
\\ \D =\frac1{-v}(\det f_{kl})^\frac1{n+2} (v+f) \ge 0.
\end{array}
\end{equation}
Here we have used (\ref{defaffinemet}), (\ref{gradv}) and Lemma \ref{positive}.

Now for the third term of (\ref{gradphi}), note that since $w$ and $\alpha = (\det f_{ij})^{-\frac1{n+2}}$ are sections of the same bundle on $M$, we can write $w = \psi \alpha$ for a positive function $\psi$ defined on $M$.
As above in Section \ref{aff-diff}, interpret $\alpha$ on $\tilde{M}$ as a component of the affine conormal vector of our hypersurface. Therefore, we have 
$$\alpha_{,11} = -\sum_k A_{11k}\alpha_{,k} - B_{11} \alpha. $$
Now compute 
$$
w_{,11}=\psi \alpha_{,11} + 2 \psi_{,1} \alpha_{,1} + \psi_{,11} \alpha \quad \mbox{and}$$
$$ \alpha = \frac{w}\psi, \quad \alpha_{,k} = \frac{\psi w_{,k} - \psi_{,k} w}{\psi^2}.$$
Since each $w_{,k}=0$ for $k\ge2$, we find from these equations that at $q$,
\begin{equation} \label{third-term}
w_{,11} \le C_1 w_{,1} + C_2 w
\end{equation}
where $C_1$ and $C_2$ depend on the norms of the cubic form $A_{ijk}$, the shape operator $B_{ij}$, and covariant derivatives of $\psi$.  Now since $M$ is compact, the affine invariants $A_{ijk}$ and $B_{ij}$, the function $\frac1{\psi}$ and the covariant derivatives of $\psi$ are all bounded in norm.  So $C_1$ and $C_2$ are uniform constants.

Putting together (\ref{gradphi}), (\ref{second-term}), (\ref{first-term}) and (\ref{third-term}), we have at $q$
$$(a^2-r^2)(h-f)^{\frac12}\left( \frac{w_{,1}}{w} \right)^2 \le 
(2a + a^2 C_1)h^{\frac12} \frac{w_{,1}}{w} + a^2C_2h^{\frac12};$$
so multiplying both sides by $(a^2-r^2)(h-f)^{\frac12}$, we find at the maximum point of $\phi$,
$$\phi^2 \le (2a+a^2C_1)h^{\frac12} \phi + a^4C_2h,$$
and therefore
$$\phi \le \sfrac12\left[(2a+a^2C_1)+\sqrt{(2a+a^2C_1)^2+4a^4C_2}\right]h^{\frac12}.$$
Now the definition of $\phi$ shows that on $\Omega_h \cap B(p,a),$
$$
\frac{\|\na w\|}{w} \le \frac{\left[(2a+a^2C_1)+\sqrt{(2a+a^2C_1)^2+4a^4C_2}\right]h^{\frac12}}{2 (h-f)^{\frac12}(a^2-r^2)},$$
and by letting $a\to\infty$, we see on all $\Omega_h$ that 
$$\frac{\|\na v\|}{-v}=\frac{\|\na w\|}{w} \le C_3(h-f)^{-\frac12},$$
where now $C_3$ depends only on $h$, $C_1$ and $C_2$.
\end{proof}

Now we use this estimate to get a contradiction. Consider the point $y$ in $\partial \Omega_h$. Consider the line segment in $x$ coordinates from the origin $p$ to $y=(y^1,\dots,y^n)$. By changing coordinates, we may assume that $y^i=0$ for $i\ge2$. Parametrize the line segment by $x(z)=zy$. 
Basically since $M$ is compact (Lemma \ref{infinite-len} below), the curve on the graph $(x(z),f(x(z)))$ has infinite length with respect to the affine metric. We will use the gradient estimate to show it must be finite.

Using the property (\ref{gradv}) of the Legendre transform, compute 
\begin{eqnarray*}
\frac{\|\na v\|^2}{v^2}&=&\frac{1}{v^2}(\det f_{kl})^{\frac1{n+2}}f^{ij}\,\frac{\partial v}{\partial x^i}\, \frac{\partial v}{\partial x^j} \\
&=&\frac{1}{v^2}(\det f_{kl})^{\frac1{n+2}}f_{ij} x^ix^j\\
&=&\psi^2(\det f_{kl})^{-\frac1{n+2}}f_{ij} x^ix^j,
\end{eqnarray*}
where $\psi = -v^{-1}(\det f_{kl})^{\frac1{n+2}}$ is a function which descends to $M$ as above.

We consider the line segment as above.  For $z \in [\frac12,1[$, consider the curve $(x(z),f(x(z)))$ as above. Its length with respect to the affine metric is 
\begin{eqnarray}
\ell&=&\int_{\frac12}^1\left[(\det f_{kl})^{-\frac1{n+2}}f_{ij} y^iy^j \right]^\frac12 dz  \label{arc-length} \\
&=&\int_{\frac12}^1 \frac1z\left[(\det f_{kl})^{-\frac1{n+2}}f_{ij} x^ix^j \right]^\frac12 dz  \nonumber \\
&=&\int_{\frac12}^1 \frac{\|\na v\|}{-zv\psi} dz \le C_4 \int_{\frac12}^1 (h-f)^{-\frac12} dz, \nonumber
\end{eqnarray}
where by the gradient estimate $C_4$ depends only on $h$, the geometry 
of $M$ and the section $u$.  Now it is a simple matter to change the 
variable of integration from $z$ to $f$ to show that this integral 
is finite. (Since $\int^h(h-f)^{-\frac12}df < \infty$.) 

\begin{lem} \label{infinite-len}
The length with respect to the affine metric of the curve $(x(z),f(x(z))),\,z \in [\frac12,1[$ is infinite.
\end{lem}

\begin{proof}
This path is a geodesic with respect to the centroaffine connection $\nabla^c$, which descends to $M$. [Just notice this ray in the $x$ coordinates is also a ray in the projective $t$ coordinates.] By the argument in Lemma \ref{maxh-lemma} above, the image on $M$ of the curve is a geodesic ray of $\nabla^c$ extended indefinitely. Since $M$ is compact and the affine metric is Riemannian, the length must be infinite.
\end{proof}
 
This lemma provides the contradiction. Therefore, each $\Omega_h \subset \subset \tilde{M}$. Consider $N=\bigcup_{h>0}\Omega_h$. Then $N$ is convex in the inhomogeneous coordinate patch $\mathcal{O}$. Since $N$ supports a negative convex function with a minimum at $p$, it must be bounded in $\mathcal{O}$. Also, $u\to 0$ at $\partial N$ implies that $\tilde{M}=N$ and Theorem \ref{main-thm} is proved.

Also, note the proof gives us the following generalization.

\begin{thm} \label{thm-boundary}
Let $M^n$ be a compact locally projectively flat manifold with boundary $\partial M$. The following statements are equivalent:
\begin{enumerate}
\item $M$ admits an oriented tautological bundle $\tau$ and there is a negative section $u$ of the the dual bundle $\tau^*$ whose Hessian matrix $u_{ij}>0$ (derivatives taken with respect to some $\rp^n$ coordinates) and $u\to0$ at the boundary of $M$.
\item $M^n$ is properly convex.
\end{enumerate}
\end{thm}

\begin{rem}
Notice that we require no regularity a priori for the boundary $\partial M$. 
\end{rem}

\begin{proof}
We only use the compactness of $M$ to get the relevant estimates.  Consider $\Omega_h$ as above. Consider the covering map $\pi:\, \tilde{M} \to M$. Then we want 
\begin{equation} \label{compact}
\pi(\Omega_h) \subset \subset M^o, \qquad \mbox{where } M^o=M\setminus \partial M.
\end{equation}
As above, we abuse the notation by denoting $u$ as the function on $\Omega_h \subset 
\tilde{M}$ developed from the section $u$ on $M$ and with a minimum 
at a specified point $p$.  Recall that $\Omega_h =\{u<-\frac1h\}$. 
Put an arbitrary metric on $M$ which preserves the topology and lift this metric to $\tilde{M}$. Then it is clear that $u$ is a continuous function on dev$^{-1}(\mathcal{O}) \subset \tilde{M}$, where $\mathcal{O}$ is the affine chart containing $\Omega_h$. 
We prove the statement (\ref{compact}) by contradiction. 
Let $q_i$ be a sequence of points in $\tilde{M}$ with $\pi(q_i) \in \partial M$ and $\lim q_i = q \in \partial \Omega_h \subset \tilde{M}$. Therefore, $u(q)=-\frac1h <0$. But $u(q)=\lim u(q_i) =0$ provides the contradiction.

Now since $\pi(\Omega_h)$ is relatively compact in $M^o$, the affine invariants $A_{ijk}$ and $B_{ij}$ and the function $\psi$ are all bounded by an amount which only depends on $h$.  Notice we still have Lemma \ref{maxh-lemma}.
Also the affine arc-length of the curve considered in Equation \ref{arc-length} above must be infinite, since we choose the path so that it cannot be extended beyond $\Omega_h$ and its image on $M$ stays outside a neighborhood of $\partial M$.  These observations complete the proof of 1 $\Rightarrow$ 2.

To show 2 $\Rightarrow$ 1, first it is easy to show that $M$ properly convex implies that $M$ admits an oriented tautological bundle (see \cite{loftin01}).
Consider $\Omega$ the bounded convex domain to which $\tilde{M}$ is projectively equivalent.
Then by Cheng-Yau \cite{cheng-yau77}, we can solve the Dirichlet problem
$$ \det (\bar{u}_{ij})=\left(-\frac1{\bar{u}}\right)^{n+2}, \qquad \bar{u}_{ij}>0, \qquad
\bar{u}|_{\partial \Omega} = 0.$$
If we treat $\bar{u}$ as a section of $\tau^*$, it descends to the quotient under any projective automorphism of $\Omega$.
\end{proof}

Note that this last argument gives this corollary:

\begin{cor}
Let $M$ be as in Theorem \ref{thm-boundary} above. Then there is a section $\bar{u}$ of $\tau^*$ on $M$ so that 
$$\det (\bar{u}_{ij})=\left(-\frac1{\bar{u}}\right)^{n+2}, \qquad \bar{u}_{ij}>0, \qquad
\bar{u}|_{\partial M} = 0.$$
Note that if $\phi = \frac{\bar{u}}u$ then $\phi$ can be thought of as perturbing $u$ to the solution to the affine sphere equation as above in Equation \ref{perturb-eq}.
\end{cor}

\bibliographystyle{abbrv}
\bibliography{thesis}

\end{document}